\documentclass[reqno,openbib,runningheads,a4paper,12pt]{amsart}%
\usepackage{graphicx}
\usepackage{amssymb}
\usepackage{amsfonts}
\usepackage{amsmath}
\usepackage{graphicx}
\usepackage{lineno}
\usepackage[authoryear]{natbib}
\usepackage[dvipsnames]{xcolor}
\usepackage{epsf}
\usepackage{lscape}
\usepackage[legalpaper,bookmarks=true,colorlinks=true,linkcolor=blue,citecolor=blue]%
{hyperref}%
\providecommand{\U}[1]{\protect\rule{.1in}{.1in}}
\providecommand{\U}[1]{\protect \rule{.1in}{.1in}}
\newtheorem{theorem}{Theorem}[section]
\newtheorem{corollary}{Corollary}[section]

\newtheorem{remark}{Remark}[section]

\makeatletter
\renewcommand{\@biblabel}[1]{}
\makeatother
\begin{document}

\begin{center}
{\Large \textbf{Nelson-Aalen kernel estimator to the tail index}}%
\newline{\Large \textbf{of right censored Pareto-type data}}\bigskip

{\large Nour Elhouda Guesmia, Abdelhakim Necir}$^{\ast}${\large , Djamel
Meraghni}\medskip\newline

{\small \textit{Laboratory of Applied Mathematics, Mohamed Khider University,
Biskra, Algeria}}\medskip\medskip
\end{center}

\noindent On the basis of Nelson-Aalen product-limit estimator of a randomly
right censored distribution function, we introduce a kernel estimator to the
tail index of right censored Pareto-like data. Under some regularity
assumptions, the consistency and asymptotic normality of the proposed
estimator are established. A simulation study shows that the proposed
estimator performs better, in terms of stability, bias and mean squared error
(MSE), than the non-smoothed one. The results are applied to insurance loss
data to illustrate their practical effectiveness.\medskip\medskip

\medskip

\noindent\textbf{AMC 2020 Subject Classification:} Primary 62G32; 62G05; 62G20.

\vfill

\vfill

\noindent{\small $^{\text{*}}$Corresponding author:
ah.necir@univ-biskra.dz\newline\noindent\textit{E-mail address:}\newline
nourelhouda.guesmia@univ-biskra.dz\texttt{\ }(N.~Guesmia)\newline
djamel.meraghni@univ-biskra.dz (D.~Meraghni)}\newline

\section{\textbf{Introduction\label{sec1}}}

\noindent Right censored data, typified by partially observed values, present
a common challenge in statistics. Indeed, this type of data are known only
within a certain range and remain unspecified outside this range. This
phenomenon is observed in fields such as actuarial science, survival analysis
and reliability engineering. Heavy-tailed distributions are used to model data
with extreme values. Compared to the normal distribution, they are
characterized by a significant probability of observing values far from the
mean. These distributions are crucial in applications like insurance, finance,
and environmental studies, where rare values can have a substantial impact.
Let $X_{1},X_{2},...,X_{n}$ be a sample of size $n\geq1$ from a random
variable (rv) $X$, and let $C_{1},C_{2},...,C_{n}$ be a sample from second rv
$C$ independent from the first one. Both rv's are defined on a probability
space $\left(  \Omega,\mathcal{A},\mathbf{P}\right)  ,$ with respective
continuous cumulative distribution functions (cdf's) $F$ and $G.$ We suppose
that $X$ is right censored by $C,$ that is, at each stage $1\leq j\leq n,$ we
can only observe the variable $Z_{j}:=\min\left(  X_{j},C_{j}\right)  $ and
the indicator variable $\delta_{j}:=\mathbb{I}_{\left\{  X_{j}\leq
C_{j}\right\}  },$ which determines whether or not $X$ has been observed. We
assume that the tail functions $\overline{F}:=1-F$ and $\overline{G}:=1-G$ are
regularly varying at infinity (or Pareto-type) of positive tail indices
$\gamma_{1}>0$ and $\gamma_{2}>0$ respectively. In other words, for any $x>0,$
we have%
\begin{equation}
\lim_{u\rightarrow\infty}\frac{\overline{F}(ux)}{\overline{F}(u)}%
=x^{-1/\gamma_{1}}\text{ and }\lim_{u\rightarrow\infty}\frac{\overline{G}%
(ux)}{\overline{G}(u)}=x^{-1/\gamma_{2}}. \label{RV}%
\end{equation}
If we denote the cdf of $Z$ by $H,$ then in virtue of the independence of $X$
and $C,$ we have $\overline{H}=\overline{F}\times\overline{G},$ which yields
that $\overline{H}$ is also regularly varying at infinity, with tail index
$\gamma:=\gamma_{1}\gamma_{2}/\left(  \gamma_{1}+\gamma_{2}\right)  .$ In the
presence of extreme values and right censored data, several methods have been
proposed for estimating the tail index. Numerous techniques have been
developed to address the challenges arising from the nature of such data. Many
studies have focused on modifying traditional tail index estimation
procedures, such as Hill's methodology \citep[][]{Hill75}, to make them
suitable for censored data. For instance, \cite{EnFG08} adapted the latter to
handle the characteristics of right censored data. Their estimator is given by
$\widehat{\gamma}_{1}^{\left(  EFG\right)  }:=\widehat{\gamma}^{\left(
H\right)  }/\widehat{p},$ where $\widehat{\gamma}^{\left(  H\right)  }:=k^{-1}%
{\displaystyle\sum\limits_{i=1}^{k}}
\log\left(  Z_{n-i+1:n}/Z_{n-k:n}\right)  ,$ denotes the popular Hill
estimator of the tail index $\gamma$ corresponding to the complete $Z$-sample
and
\begin{equation}
\widehat{p}:=k^{-1}%
{\displaystyle\sum_{i=1}^{k}}
\delta_{\left[  n-i+1:n\right]  }, \label{p=hate}%
\end{equation}
stands for an estimator of the proportion of upper non-censored observations%
\begin{equation}
p:=\gamma_{2}/\left(  \gamma_{1}+\gamma_{2}\right)  . \label{p}%
\end{equation}
The integer sequence $k=k_{n}$ represents the number of top order statistics
satisfying
\begin{equation}
k\rightarrow\infty\text{ and }k/n\rightarrow0\text{ as }n\rightarrow\infty.
\label{k}%
\end{equation}
The sequence of rv's $Z_{1:n}\leq...\leq Z_{n:n}$ represent the order
statistics pertaining to the sample $Z_{1},...,Z_{n}$ and $\delta_{\left[
1:n\right]  },...,\delta_{\left[  n:n\right]  }$ denote the corresponding
concomitant values, satisfying $\delta_{\left[  j:n\right]  }=\delta_{i}$ for
$i$ such that $Z_{j:n}=Z_{i}.$ Useful Gaussian approximations to
$\widehat{p},$ $\widehat{\gamma}^{\left(  H\right)  }$ and $\widehat{\gamma
}_{1}^{\left(  EFG\right)  }$ in terms of sequences of Brownian bridges are
given in \cite{BMN-2015}. On the basis of Kaplan-Meier integration,
\cite{WW2014}, proposed a consistent estimator to $\gamma_{1}$ defined by
\[
\widehat{\gamma}_{1}^{\left(  W\right)  }:=%
{\displaystyle\sum\limits_{i=1}^{k}}
\frac{\overline{F}_{n}^{\left(  KM\right)  }\left(  Z_{n-i:n}\right)
}{\overline{F}_{n}^{\left(  KM\right)  }\left(  Z_{n-k:n}\right)  }\log
\frac{Z_{n-i+1:n}}{Z_{n-i:n}},
\]
where $F_{n}^{\left(  KM\right)  }\left(  z\right)  :=1-%
{\displaystyle\prod_{Z_{i:n}\leq z}}
\left(  \dfrac{n-i}{n-i+1}\right)  ^{\delta_{\left[  i:n\right]  }}$ denotes
the popular Kaplan-Meier estimator of cdf $F$ \citep[] []{KM58}. The
simulation study given by \cite{BMV2018} showed that $\widehat{\gamma}%
_{1}^{\left(  W\right)  }$ has the best MSE behaviour between all available
first order estimators of $\gamma_{1}.$ The asymptotic normality of the latter
is established in \cite{BWW2019} by considering Hall's model
\citep[] []{Hall82}. The bias reduction in tail index estimation for right
censored data is addressed by \cite{BBWG2016} and \cite{BWW2019}. By using a
Nelson-Aalen integration, recently \cite{MNS2025} also derived a new
$\gamma_{1}$ estimator%
\[
\widehat{\gamma}_{1}^{\left(  MNS\right)  }:=%
{\displaystyle\sum\limits_{i=1}^{k}}
\frac{\delta_{\left[  n-i+1:n\right]  }}{i}\frac{\overline{F}_{n}^{\left(
NA\right)  }\left(  Z_{n-i+1:n}\right)  }{\overline{F}_{n}^{\left(  NA\right)
}\left(  Z_{n-k:n}\right)  }\log\frac{Z_{n-i+1:n}}{Z_{n-k:n}},
\]
where
\begin{equation}
F_{n}^{\left(  NA\right)  }\left(  z\right)  =1-%
{\displaystyle\prod_{Z_{i:n}\leq z}}
\exp\left\{  -\frac{\delta_{\left[  i:n\right]  }}{n-i+1}\right\}  ,
\label{Nelson}%
\end{equation}
denotes the well-known Nelson-Aalen estimator of cdf $F$
\citep[] []{Nelson1972}. The authors showed that both tail index estimators
$\widehat{\gamma}_{1}^{\left(  W\right)  }$ and $\widehat{\gamma}_{1}^{\left(
MNS\right)  }$ exhibit similar performances with respect to bias and mean
squared error. Actually, comparison studies made between Kaplan-Meier and
Nelson-Aalen estimators reached the conclusion that they exhibit almost
similar statistical behaviors, see, for instance, \cite{Colosimo2002}. On the
other hand, reaching to establish the asymptotic properties of extreme
Kaplan-Meier integrals poses some difficulties. To overcome this issue,
\cite{MNS2025} introduced Nelson-Aalen tail-product limit process and
established its Gaussian approximation in terms of standard Wiener processes.
Thus, they established the consistency and asymptotic normality of their
estimator $\widehat{\gamma}_{1}^{\left(  MNS\right)  }$ by considering both
regular variation conditions $\left(  \ref{RV}\right)  $ and $\left(
\ref{second-order}\right)  .$ They also showed that (as expected\textbf{)
}both the asymptotic bias and variance of $\widehat{\gamma}_{1}^{\left(
W\right)  }$ and $\widehat{\gamma}_{1}^{\left(  MNS\right)  }$ are
respectively equal.\medskip

\noindent It is well known that, in the usual complete data case, introducing
a kernel function in classical statistical estimation produces estimators with
nice properties in terms of stability and bias. This is also true in extreme
value theory (EVT) based estimation where, for instance, \cite{CDM85} and
\cite{Wolf 2003} proposed smooth tail index estimators. In the case of
truncated heavy-tailed data, this problem was addressed by
\cite{Benchaira2016b} who proposed a smooth version to the estimator they
defined in \cite{Benchaira2016a}. To the best of our knowledge, kernel
functions have not been used yet to estimate the extreme value index of
censored datasets. Our main goal in this work is precisely to deal with this
situation by providing a kernel version to $\widehat{\gamma}_{1}^{\left(
MNS\right)  }$ and establishing its consistency and asymptotic normality
thanks to the aforementioned Gaussian approximation. To this end, let us
define a specific function $K:\mathbb{R}\rightarrow\mathbb{R}_{+}$ called the
kernel, which satisfies the following conditions:

\begin{itemize}
\item $\left[  A1\right]  $ is non-increasing and right-continuous on
$\mathbb{R}.$

\item $\left[  A2\right]  $ $K\left(  s\right)  =0$ for $s\notin\left[
0,1\right)  $ and $K\left(  s\right)  \geq0$ for $s\in\left[  0,1\right)  .$

\item $\left[  A3\right]  $ $\int_{\mathbb{R}}K\left(  s\right)  ds=1.$

\item $\left[  A4\right]  $ $K$ and its first and second Lebesgue derivatives
$K^{\prime}$ and $K^{\prime\prime}$ are bounded on $\mathbb{R}.$\medskip
\end{itemize}

\noindent Such functions include, for example, the indicator kernel
$K_{1}:=\mathbb{I}_{\left\{  0\leq s<1\right\}  },$ the biweight $K_{2}$ and
triweight kernels $K_{3},$ which are defined by
\begin{equation}
K_{2}\left(  s\right)  :=\frac{15}{8}\left(  1-s^{2}\right)  ^{2}\text{ and
}K_{3}\left(  s\right)  :=\frac{35}{16}\left(  1-s^{2}\right)  ^{3},\text{
}0\leq s<1, \label{K2-K3}%
\end{equation}
and zero elsewhere, see, e.g., \cite{Wolf 2003}. By applying Potter's
inequalities, see e.g. Proposition B.1.10 in \cite{deHF06}, to the regularly
varying function $F$ together with assumptions $[A1]-[A2],$
\cite{Benchaira2016b} stated that%
\begin{equation}
\lim_{u\rightarrow\infty}\int_{u}^{\infty}w_{K}\left(  \frac{\overline
{F}\left(  y\right)  }{\overline{F}\left(  u\right)  }\right)  \log\frac{y}%
{u}d\frac{F\left(  y\right)  }{\overline{F}\left(  u\right)  }=\gamma_{1}%
\int_{0}^{\infty}K\left(  y\right)  dy=\gamma_{1}, \label{gamma1}%
\end{equation}
where $w_{K}\left(  s\right)  :=d\left\{  sK\left(  s\right)  \right\}  /ds$
denotes the Lebesgue derivative of function $s\rightarrow sK\left(  s\right)
.$ By letting $u=Z_{n-k:n}$ and substituting $F$ by Nelson-Aalen estimator
$F_{n}^{\left(  NA\right)  },$ we derive a kernel estimator to the tail index
$\gamma_{1}$ as follows%
\begin{equation}
\widehat{\gamma}_{1,K}:=\int_{Z_{n-k:n}}^{\infty}w_{K}\left(  \frac
{\overline{F}_{n}^{\left(  NA\right)  }\left(  y\right)  }{\overline{F}%
_{n}^{\left(  NA\right)  }\left(  Z_{n-k:n}\right)  }\right)  \log\frac
{y}{Z_{n-k:n}}d\frac{F_{n}^{\left(  NA\right)  }\left(  y\right)  }%
{\overline{F}_{n}^{\left(  NA\right)  }\left(  Z_{n-k:n}\right)  }.
\label{f-form}%
\end{equation}
Let us rewrite the previous integral into a more friendly form as a sum of
terms from the $Z$-sample. To this end, we use the following crucial
representation of cdf $F$ in terms of estimable functions $H$ and $H^{\left(
1\right)  }:$%
\[
\int_{0}^{z}\frac{dH^{\left(  1\right)  }\left(  x\right)  }{\overline
{H}\left(  x-\right)  }=\int_{0}^{z}\frac{dF\left(  x\right)  }{\overline
{F}\left(  x\right)  }=-\log\overline{F}\left(  x\right)  ,\text{ for }z>0,
\]
where $H^{\left(  1\right)  }\left(  z\right)  :=\mathbf{P}\left(  Z\leq
z,\delta=1\right)  =\int_{0}^{z}\overline{G}\left(  x\right)  dF\left(
x\right)  .$ This implies that%
\begin{equation}
\overline{F}\left(  z\right)  =\exp\left\{  -\int_{0}^{z}\frac{dH^{\left(
1\right)  }\left(  x\right)  }{\overline{H}\left(  x-\right)  }\right\}
,\text{ for }z>0, \label{F-formula}%
\end{equation}
where $\overline{H}\left(  x-\right)  :=\lim_{\epsilon\downarrow0}$
$\overline{H}\left(  x-\epsilon\right)  .$ The empirical counterparts of cdf
$H$ and sub-distribution $H^{\left(  1\right)  }$ are given by
\[
H_{n}\left(  z\right)  :=n^{-1}\sum_{i=1}^{n}\mathbb{I}_{\left\{  Z_{i:n}\leq
z\right\}  }\text{ and }H_{n}^{\left(  1\right)  }\left(  z\right)
:=n^{-1}\sum_{i=1}^{n}\delta_{\left[  i:n\right]  }\mathbb{I}_{\left\{
Z_{i:n}\leq z\right\}  },
\]
see, for instance, \cite{SW86} page 294. Note that $\overline{H}_{n}\left(
Z_{n:n}\right)  =0,$ for this reason and to avoid dividing by zero, we use
$\overline{H}\left(  x-\right)  $ instead of $\overline{H}\left(  x\right)  $
in formula $\left(  \ref{F-formula}\right)  .$ Nelson-Aalen estimator
$F_{n}^{\left(  NA\right)  }\left(  z\right)  $ is obtained by substituting
both $H_{n}$ and $H_{n}^{\left(  1\right)  }$ in the right-hand side of
formula $\left(  \ref{F-formula}\right)  .$ That is, we have
\begin{equation}
\overline{F}_{n}^{\left(  NA\right)  }\left(  z\right)  :=\exp\left\{
-\int_{0}^{z}\frac{dH_{n}^{\left(  1\right)  }\left(  x\right)  }{\overline
{H}_{n}\left(  x-\right)  }\right\}  . \label{emp-F-formula}%
\end{equation}
Since $\overline{H}_{n}\left(  Z_{i:n}-\right)  =\overline{H}_{n}\left(
Z_{i-1:n}\right)  =\left(  n-i+1\right)  /n,$ then the previous integral leads
to the expression of $\overline{F}_{n}^{\left(  NA\right)  }\left(  z\right)
$ in $\left(  \ref{Nelson}\right)  .$ Differentiating both sides of formula
$\left(  \ref{emp-F-formula}\right)  $ yields that%
\[
dF_{n}^{\left(  NA\right)  }\left(  z\right)  =\frac{\overline{F}_{n}^{\left(
NA\right)  }\left(  z\right)  }{\overline{H}_{n}\left(  z-\right)  }%
dH_{n}^{\left(  1\right)  }\left(  z\right)  .
\]
By substituting $dF_{n}^{\left(  NA\right)  }\left(  z\right)  $ by the above
in $\left(  \ref{f-form}\right)  ,$ we end up with%
\begin{equation}
\widehat{\gamma}_{1,K}=\sum_{i=1}^{k}\frac{\delta_{\left[  n-i+1:n\right]  }%
}{i}\frac{\overline{F}_{n}^{\left(  NA\right)  }\left(  Z_{n-i+1:n}\right)
}{\overline{F}_{n}^{\left(  NA\right)  }\left(  Z_{n-k:n}\right)  }%
w_{K}\left(  \frac{\overline{F}_{n}^{\left(  NA\right)  }\left(
Z_{n-i+1:n}\right)  }{\overline{F}_{n}^{\left(  NA\right)  }\left(
Z_{n-k:n}\right)  }\right)  \log\frac{Z_{n-i+1:n}}{Z_{n-k:n}}.
\label{gamma1-estim}%
\end{equation}
The remainder of the paper is structured as follows. In Section \ref{sec2}, we
present our main results namely the consistency and asymptotic normality of
$\widehat{\gamma}_{1,K},$ whose proofs are postponed to Section \ref{sec5}. We
devote Section \ref{kopt} to the crucial issue of choosing the optimal number
$k$ of upper order statistics used in estimate computation. In Section
\ref{sec3}, we carry out a simulation study to illustrate the finite sample
behavior of our estimator with a comparison with the non-smoothed ones
$\widehat{\gamma}_{1}^{\left(  MNS\right)  }.$ In Section \ref{sec4}, we apply
our results to a real dataset of insurance losses.

\section{\textbf{Main results\label{sec2}}}

\noindent Since weak approximations of extreme value theory based statistics
are achieved in the second-order framework
\citep[see e.g. page 48  in] []{deHF06}, then it seems quite natural to
suppose that cdf $F$ satisfies the well-known second-order condition of
regular variation. That is, we assume that for any $x>0:$%
\begin{equation}
\underset{t\rightarrow\infty}{\lim}\frac{U_{F}\left(  tx\right)  /U_{F}\left(
t\right)  -x^{\gamma_{1}}}{A_{1}^{\ast}\left(  t\right)  }=x^{\gamma_{1}%
}\dfrac{x^{\tau_{1}}-1}{\tau_{1}}, \label{second-order}%
\end{equation}
or equivalently%
\begin{equation}
\lim_{t\rightarrow\infty}\frac{\overline{F}\left(  tx\right)  /\overline
{F}\left(  t\right)  -x^{-1/\gamma_{1}}}{A_{1}\left(  t\right)  }%
=x^{-1/\gamma_{1}}\frac{x^{\tau_{1}/\gamma_{1}}-1}{\tau_{1}\gamma_{1}},
\label{second-orderF}%
\end{equation}
where $\tau_{1}\leq0$\ is the second-order parameter,and $A_{1}\left(
t\right)  :=A_{1}^{\ast}\left(  1/\overline{F}\left(  t\right)  \right)  $
with $A_{1}^{\ast}$ being a function tending to $0,$ not changing sign near
infinity and having a regularly varying absolute value at infinity with index
$\tau_{1}.$\ If $\tau_{1}=0,$ interpret $\left(  x^{\tau_{1}}-1\right)
/\tau_{1}$ as $\log x.$ The functions $L^{\leftarrow}\left(  s\right)
:=\inf\left\{  x:L\left(  x\right)  \geq s\right\}  ,$ for $0<s<1,$ and
$U_{L}\left(  t\right)  :=L^{\leftarrow}\left(  1-1/t\right)  ,$ $t>1,$
respectively stand for the quantile and tail quantile functions of any cdf
$L.$ A very popular class of such distributions consists of what is called
Hall's model \citep[][]{Hall82}. It includes the most usual Pareto-type cdf's
namely, the distributions of Burr and Fr\'{e}chet, the generalized extreme
value distribution (GEV), the generalized Pareto distribution (GPD), etc. In
real-life applications, these probability distributions proved to be extremely
useful for risk managers as they provide very appropriate statistical models
for large losses in various fields, such as insurance, finance, hydrology and
environmental sciences. Hall's model is characterized by its tail%
\begin{equation}
\overline{F}\left(  x\right)  =C_{1}x^{-1/\gamma_{1}}\left(  1+C_{2}%
x^{\tau_{1}/\gamma_{1}}+o\left(  x^{\tau_{1}/\gamma_{1}}\right)  \right)
,\text{ as }x\rightarrow\infty, \label{HallF}%
\end{equation}
where $C_{1}\neq0,$ $C_{2}\in\mathbb{R}$ and $\tau_{1}<0.$ We show that
$\overline{F}$ satisfies the second order condition $\left(
\ref{second-orderF}\right)  $ with
\begin{equation}
A_{1}\left(  t\right)  =\left(  1+o\left(  1\right)  \right)  \left(  \tau
_{1}/\gamma_{1}\right)  C_{2}t^{\tau_{1}/\gamma_{1}}, \label{A1}%
\end{equation}
see for instance page 77 in the textbook of \cite{deHF06}. For example:

\begin{itemize}
\item Burr's model $\overline{F}\left(  x\right)  =\left(  1+x^{1/\eta
}\right)  ^{-\eta/\gamma_{1}}$ for $\eta,\gamma_{1}>0:$%
\[
\overline{F}\left(  x\right)  =x^{-1/\gamma_{1}}\left(  1-\frac{\eta}%
{\gamma_{1}}x^{-1/\eta}+o\left(  x^{-1/\eta}\right)  \right)  ,\text{ as
}x\rightarrow\infty.
\]
with $\tau_{1}=-\gamma_{1}/\eta,$ $C_{1}=1,$ $C_{2}=-\eta/\gamma_{1}$ and
$A_{1}\left(  t\right)  :=\gamma_{1}^{-1}\left(  1+o\left(  1\right)  \right)
t^{-1/\eta}.$

\item Fr\'{e}chet model $\overline{F}\left(  x\right)  =1-\exp\left(
-x^{-\alpha}\right)  ,$ for $\alpha>0:$%
\[
\overline{F}\left(  x\right)  =x^{-\alpha}\left(  1-\frac{1}{2}x^{-\alpha
}+o\left(  x^{-\alpha}\right)  \right)  .
\]
with $\gamma_{1}=\alpha^{-1},$ $\tau_{1}=-1,$ $C_{1}=1,$ $C_{2}=-1/2$ and
$A_{1}\left(  t\right)  =\left(  1+o\left(  1\right)  \right)  \left(
\alpha/2\right)  t^{-\alpha}.$
\end{itemize}

\noindent Let us now state our main results.

\begin{theorem}
\textbf{\label{theorem}}Assume that cdf's $F$ and $G$ satisfy the first-order
condition of regular variation $\left(  \ref{RV}\right)  $ with $\gamma
_{1}<\gamma_{2},$ and let $K$ be a kernel function satisfying assumptions
$[A1]-[A4].$ For an integer sequence $k=k_{n}$ such that $k\rightarrow\infty$
and $k/n\rightarrow0,$ we have $\widehat{\gamma}_{1,K}\overset{\mathbf{P}%
}{\rightarrow}\gamma_{1},$ as $n\rightarrow\infty.$ In addition, if we set
$h:=U_{H}\left(  n/k\right)  $ and assume that the second-order condition of
regular variation $\left(  \ref{second-order}\right)  $ holds, then%
\begin{equation}
\sqrt{k}\left(  \widehat{\gamma}_{1,K}-\gamma_{1}\right)  \overset{\mathcal{D}%
}{=}\mathcal{N}\left(  0,\sigma_{K}^{2}\right)  +\sqrt{k}A_{1}\left(
h\right)  \mu_{K}+o_{\mathbf{P}}\left(  1\right)  ,\text{ as }n\rightarrow
\infty, \label{approx1}%
\end{equation}
provided that $\sqrt{k}A_{1}\left(  h\right)  =O\left(  1\right)  ,$ where%
\[
\mu_{K}:=\int_{0}^{1}s^{-\tau_{1}}K\left(  s\right)  ds\text{ and }\sigma
_{K}^{2}:=\gamma_{1}^{2}\int_{0}^{1}s^{-1/p+1}K^{2}\left(  s\right)  ds.
\]
In addition, If we assume that $\sqrt{k}A_{1}\left(  h\right)  \rightarrow
\lambda<\infty,$ then%
\[
\sqrt{k}\left(  \widehat{\gamma}_{1,K}-\gamma_{1}\right)  \overset{\mathcal{D}%
}{\rightarrow}\mathcal{N}\left(  \lambda\mu_{K},\sigma_{K}^{2}\right)  ,\text{
as }n\rightarrow\infty.
\]

\end{theorem}

\begin{remark}
It is clear that $\mu_{K_{1}}=1/\left(  1-\tau_{1}\right)  $ and
$\sigma_{K_{1}}^{2}=p\gamma_{1}^{2}/\left(  2p-1\right)  $ respectively
coincide with the asymptotic mean and variance of $\widehat{\gamma}%
_{1}^{\left(  MNS\right)  }$ and $\widehat{\gamma}_{1}^{\left(  W\right)  }.$
\end{remark}

\begin{remark}
In the case of complete data, that is when $p=1,$ $\lambda\mu_{K}$ and
$\sigma_{K}^{2}$ respectively coincide with the asymptotic mean and variance
of CDM's estimator \citep[] []{CDM85} defined by
\[
\widehat{\gamma}_{K}^{\left(  CDM\right)  }:=\sum_{i=1}^{k}w_{K}\left(
\frac{i}{k}\right)  \log\left(  Z_{n-i+1:n}/Z_{n-k:n}\right)  .
\]

\end{remark}

\begin{remark}
It is worth remembering that kernel estimation procedure reduces the
asymptotic bias. However it increases the asymptotic variance, which is, in
general, the price de pay when reducing the bias.
\end{remark}

\begin{remark}
\textbf{\label{remark4}}When considering the product limit estimators
$F_{n}^{\left(  KM\right)  }$ or $F_{n}^{\left(  NA\right)  }$ in the process
of estimating the tail index of Pareto-type models, the condition $\gamma
_{1}<\gamma_{2}$ arises, see for instance \cite{WW2014}, \cite{BMV2018},
\cite{BWW2019}, \cite{GM2024} and \cite{Bladt2024}. Equivalent to $p>1/2,$
this condition means that the proportion of upper non-censored observations
has to be greater than $50\%$ which seems to be quite reasonable.
\end{remark}

\section{\textbf{Optimal sample fraction\label{kopt}}}

\noindent It is clear that EVT-based estimators of the tail index and related
quantities essentially rely on the number $k$ of upper order statistics,
satisfying $\left(  \ref{k}\right)  ,$ and so their behaviors are affected by
this crucial number. Indeed, using too many data, in the estimation procedure,
results in a substantial bias whereas using too few observations leads to a
considerable variance. Then, it is legitimate to ask how to balance between
these two vices and therefore select the optimal $k$-value, called optimal
sample fraction and denoted by $k_{opt},$ in order to guarantee the best
possible estimates. In other words, the number $k_{opt}$ locates where the
distribution tail (really) begins. There are several approaches to answer this
question. A first one is based on graphical considerations. We plot the
estimator as function of $k$ and then we look for a stable region of the graph
near the horizontal line representing the true value of the parameter. That
region corresponds to the desired optimal number of top statistics. A second
method is based upon the idea of minimizing the asymptotic mean squared error
(AMSE) which is a combination of bais and variance terms. More precisely, as
function of $k,$ the AMSE is equal to the squared asymptotic bias plus the
asymptotic variance. The $k$-value where the minimum is reached represents an
optimal choice for the number of order statistics that should be used in the
estimation. Unfortunately, this process is made difficult by the fact that
this number does not depend exclusively on the sample size and the parameter
but it also depends on some unknown new parameter, namely the second order
parameter This led to several research works proposting numerical adaptive
procedures of choosing the right value of $k,$ among which we can cite, for
instance, \cite{Hall90}, \cite{DH93}, \cite{RT97} and \cite{ChP01}. In this
section, we exploit the result of Theorem \ref{theorem} to deduce, in
Corollary \ref{theorembis}, the optimal sample fraction needed for the
computation of \ the tail index estimate value $\left(  \ref{gamma1-estim}%
\right)  $ when the censoring distribution and the one which is censored both
belong to Hall's model. To this end, we also assume that, cdf $G$ satisfy
Hall's condition%
\begin{equation}
\overline{G}\left(  x\right)  =D_{1}x^{-1/\gamma_{2}}\left(  1+D_{2}%
x^{\tau_{2}/\gamma_{2}}+o\left(  x^{\tau_{2}/\gamma_{2}}\right)  \right)
,\text{ as }x\rightarrow\infty, \label{HallG}%
\end{equation}
where $D_{1}\neq0,$ $D_{2}\in\mathbb{R}$ and $\tau_{2}<0.$ It is noteworthy
that the theoretical formula of $k_{opt},$ in the corollary, is mainly useful
for illustrative simulations to check the finite sample behavior of the
estimator, as it will be seen in the next section.

\begin{corollary}
\label{theorembis}Assume that cdf's $F$ and $G$ satisfy Hall's conditions
$\left(  \ref{HallF}\right)  $ and $\left(  \ref{HallG}\right)  .$ Then, the
optimal sampe fraction is%
\[
k_{opt}:=\left[  \omega n^{-2p\tau_{1}/\left(  1-2p\tau_{1}\right)  }\right]
,
\]
where%
\[
\omega:=\left\{  -\frac{1}{2\tau_{1}^{3}}\frac{\gamma_{1}^{4}}{C_{2}%
^{2}\left(  D_{1}C_{1}\right)  ^{2p\tau_{1}}}\left(  \int_{0}^{1}%
s^{-1/p+1}K^{2}\left(  s\right)  ds\right)  \left(  \int_{0}^{1}s^{-\tau_{1}%
}K\left(  s\right)  ds\right)  ^{-2}\right\}  ^{1/\left(  1-2p\tau_{1}\right)
},
\]
with $\left[  x\right]  $ standing for the integer part of the real number
$x.$
\end{corollary}

\noindent As an application example of Corollary \ref{theorembis}, we
summarize in Table \ref{cor-kopt} some computations of $k_{opt}$ with samples
of size $1000.$%

\begin{table}[h] \centering
\begin{tabular}
[c]{ccccccccc}%
\multicolumn{5}{c}{Burr-Burr} & \multicolumn{4}{|c}{Fr\'{e}chet-Fr\'{e}chet}%
\\\hline
\multicolumn{3}{c}{$K_{1}\left(  s\right)  $} & \multicolumn{2}{|c}{$K_{2}%
\left(  s\right)  $} & \multicolumn{2}{|c}{$K_{1}\left(  s\right)  $} &
\multicolumn{2}{|c}{$K_{2}\left(  s\right)  $}\\\hline
$p$ & \multicolumn{1}{||c}{$0.6$} & \multicolumn{1}{|c}{$0.9$} &
\multicolumn{1}{|c}{$0.6$} & \multicolumn{1}{|c}{$0.9$} &
\multicolumn{1}{|c}{$0.6$} & \multicolumn{1}{|c}{$0.9$} &
\multicolumn{1}{|c}{$0.6$} & $0.9$\\\hline
\multicolumn{9}{c}{$\gamma_{1}=0.4$}\\\hline
$k_{opt}$ & \multicolumn{1}{||c}{$50$} & \multicolumn{1}{|c}{$82$} &
\multicolumn{1}{|c}{$112$} & \multicolumn{1}{|c}{$132$} &
\multicolumn{1}{|c}{$35$} & \multicolumn{1}{|c}{$50$} &
\multicolumn{1}{|c}{$81$} & \multicolumn{1}{|c}{$82$}\\\hline
\multicolumn{9}{c}{$\gamma_{1}=0.7$}\\\hline
$k_{opt}$ & \multicolumn{1}{||c}{$236$} & \multicolumn{1}{|c}{$300$} &
\multicolumn{1}{|c}{$485$} & \multicolumn{1}{|c}{$466$} &
\multicolumn{1}{|c}{$96$} & \multicolumn{1}{|c}{$112$} &
\multicolumn{1}{|c}{$224$} & \multicolumn{1}{|c}{$183$}\\\hline
&  &  &  &  &  &  &  &
\end{tabular}
\caption{Optimal sample fractions computed with samples of size 1000 under two censoring scenarios}\label{cor-kopt}%
\end{table}%

\section{\textbf{Simulation study\label{sec3}}}

\noindent To evaluate the performance of our estimator $\widehat{\gamma}%
_{1,K}$ and compare it with $\widehat{\gamma}_{1}^{\left(  MNS\right)  },$ we
consider two scenarios:

\begin{itemize}
\item Fr\'{e}chet model with parameter $\gamma_{1}$ censored by Fr\'{e}chet
model with parameter $\gamma_{2}.$

\item Burr model with parameters $(\gamma_{1},\eta)$ censored by Burr model
with parameters $(\gamma_{2},\eta).$
\end{itemize}

\begin{figure}[h]%
\centering
\includegraphics[
height=2.6844in,
width=5.482in
]%
{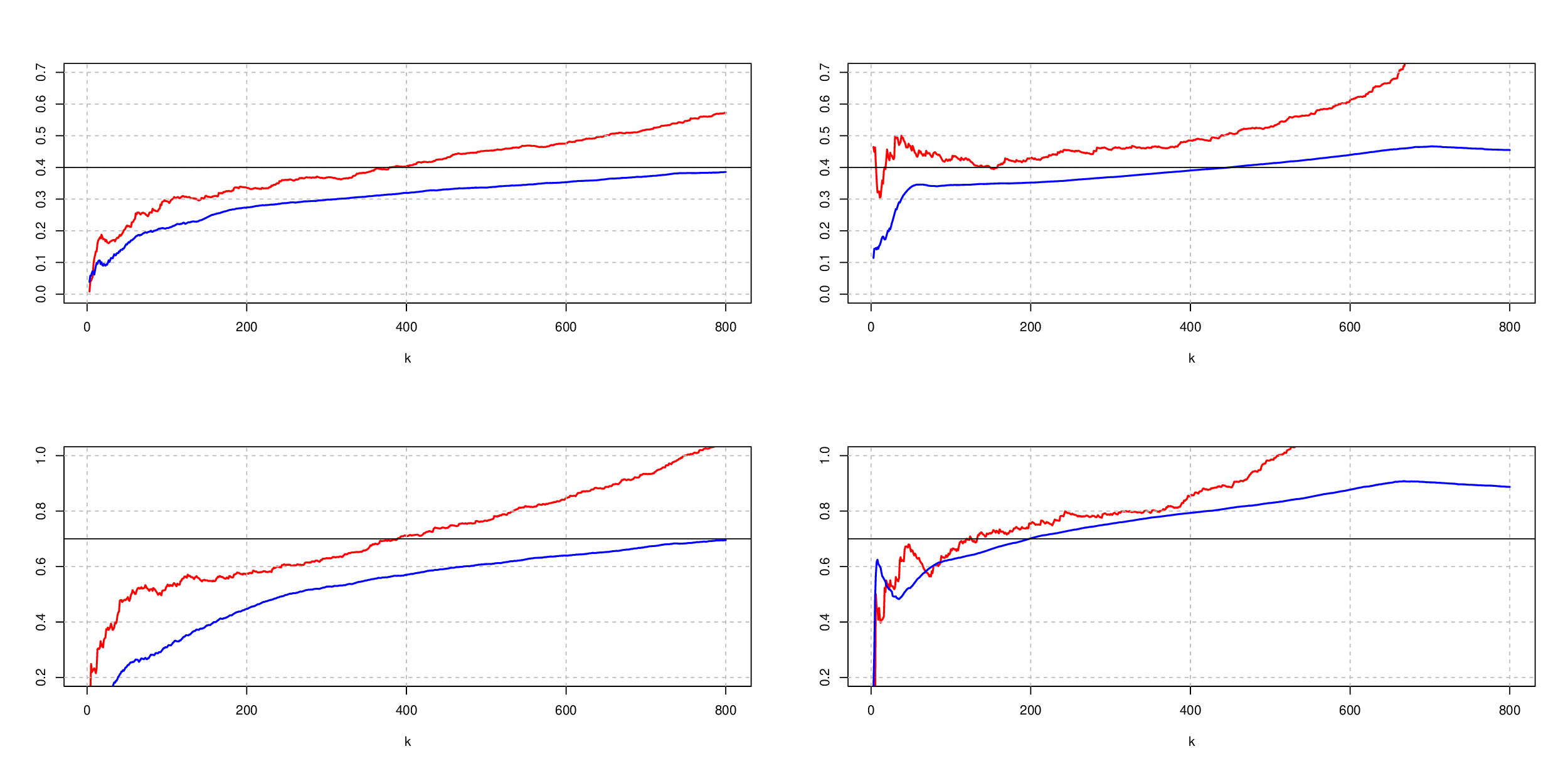}%
\caption{Plots of $\protect\widehat{\gamma}_{1,K_{2}}$ (blue line) and
$\protect\widehat{\gamma}_{1}^{\left(  MNS\right)  }$ (red line), as functions
of the number $k$ of upper order statistics, based on $2000$ samples of size
$1000$ from Fr\'{e}chet model censored by Fr\'{e}chet for $\gamma_{1}=0.4$
(top) and $\gamma_{1}=0.7$ (bottom) with $p=0.60$ (left panel) and $p=0.90$
(right panel). The horizontal line represents the true value of $\gamma_{1}.$}%
\label{fig1}%
\end{figure}

\begin{figure}[h]%
\centering
\includegraphics[
height=2.6844in,
width=5.482in
]%
{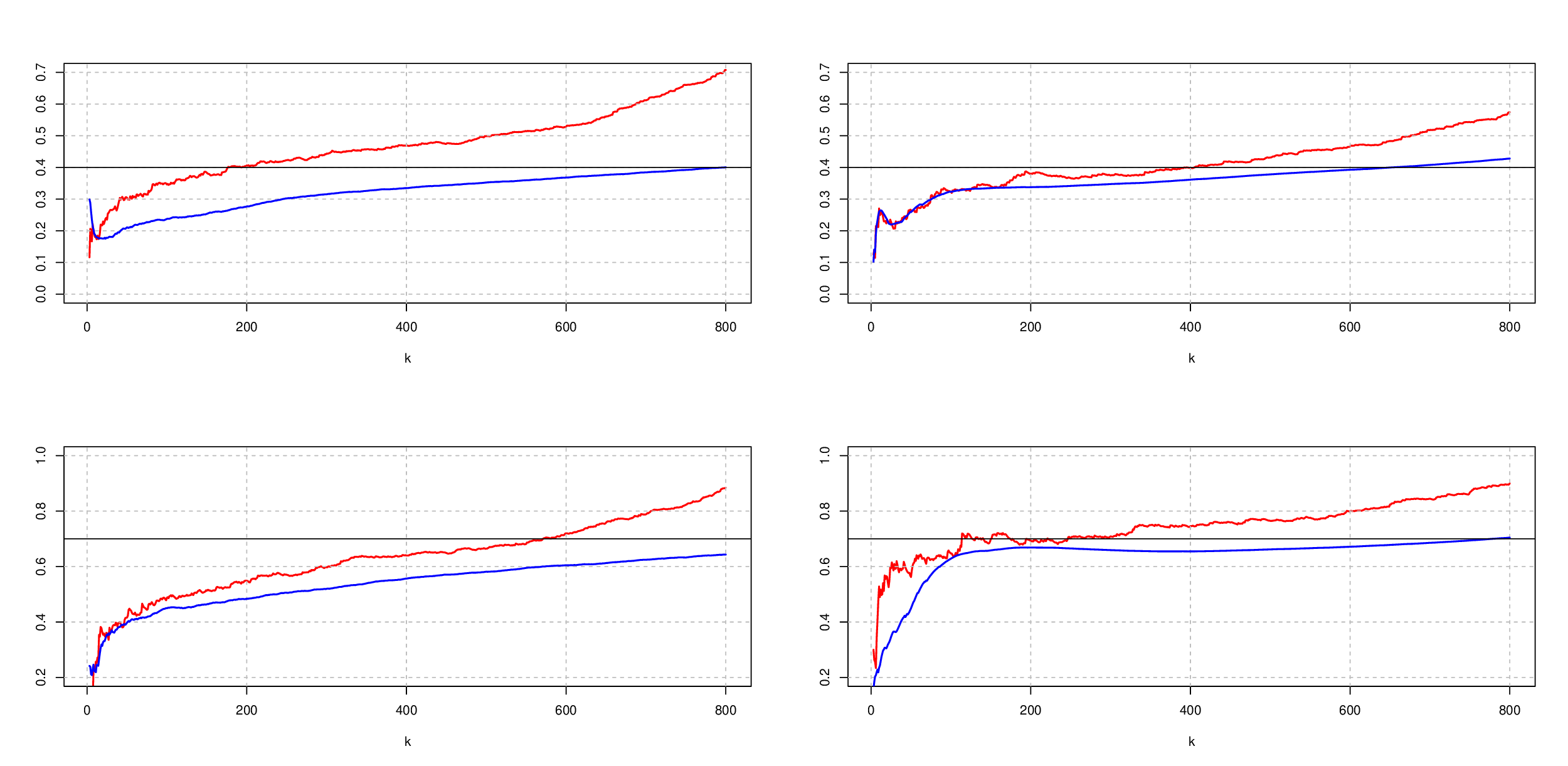}%
\caption{Plots of $\protect\widehat{\gamma}_{1,K_{2}}$ (blue line) and
$\protect\widehat{\gamma}_{1}^{\left(  MNS\right)  }$ (red line), as functions
of the number $k$ of upper order statistics, based on $2000$ samples of size
$1000$ from Burr model censored by Burr for $\gamma_{1}=0.4$ (top) and
$\gamma_{1}=0.7$ (bottom) with $p=0.60$ (left panel) and $p=0.90$ (right
panel). The horizontal line represents the true value of $\gamma_{1}.$}%
\label{fig2}%
\end{figure}

\noindent We fix $\eta=$ $0.25,\ $consider two values $0.4$ and $0.7$ for
$\gamma_{1}$ and choose two proportion values $p=0.60$ and $0.90.$ For each
combination of $\gamma_{1}$ and $p,$ we solve the equation $p=\gamma
_{2}/\left(  \gamma_{1}+\gamma_{2}\right)  $ to get the corresponding value of
$\gamma_{2}.$ We construct our estimator $\widehat{\gamma}_{1,K_{2}}$ which is
a version of $\widehat{\gamma}_{1,K}$ corresponding to the biweight kernel
function $K_{2},$ defined on the left-hand side in $\left(  \ref{K2-K3}%
\right)  $. We generate $2000$ samples of size $1000$ and our results,
obtained by averaging over all replications, are presented in Figures
\ref{fig1} and \ref{fig2} and summarized in Tables \ref{BB} and \ref{FF}. In
these the numbers $k_{opt}$ are determined by the method of Reiss and Thomas
(see \cite{RT97}, page 121) who define the optimal number of extremes needed
in estimate computation, by%
\[
k_{opt}:=\arg\min_{1<k<n}\frac{1}{k}\sum_{i=1}^{k}i^{\theta}\left\vert
\widehat{\gamma}_{1,i}-\text{median}\left\{  \widehat{\gamma}_{1,1}%
,...,\widehat{\gamma}_{1,k}\right\}  \right\vert ,\text{ }0\leq\theta\leq1/2,
\]
where $\widehat{\gamma}_{1,i}$ is $\gamma_{1}$-estimator $(\widehat{\gamma
}_{1,K_{2}}$ or $\widehat{\gamma}_{1}^{\left(  MNS\right)  })$ based on the
$i$ upper order statistics. In this study, we choose $\theta=0.3$ which
provides better results in terms of bias and MSE.\smallskip

\noindent As expected, we note that, unlike $\widehat{\gamma}_{1}^{\left(
MNS\right)  },$ our estimator enjoys smoothness and bias reduction. Moreover,
we can see in all four graphs of each figure that, compared to
$\widehat{\gamma}_{1}^{\left(  MNS\right)  },$ our estimator proves superior
stability through the range of $k.$ This is a very good and conveniently
suitable asset for people dealing with real-world data in case studies, where
only single datasets are available and no replications are possible as in
simulation studies. Finally, the top-left graphs in both Figures \ref{fig1}
and \ref{fig2} show that the accuracy of our estimator improves as the
censoring percentage and the tail index value decrease. On the other hand, the
results of Tables \ref{BB} and \ref{FF} confirm the (expected) edge that our
kernel estimator $\widehat{\gamma}_{1,K_{2}}$ has over $\widehat{\gamma}%
_{1}^{\left(  MNS\right)  }$ with respect to bias and MSE. Moreover, during
the simulation stages, we noted that the individual optimal sample fractions,
relative to $\widehat{\gamma}_{1,K_{2}},$ obtained at each one of the one
hundred replications have no significant difference unlike those pertaining to
$\widehat{\gamma}_{1}^{\left(  MNS\right)  }.$ This asset is justified by the
stability quality of our estimator.%

%

%

\begin{table}[h] \centering
\begin{tabular}
[c]{|c|c|c|c|c|}\hline
\multicolumn{5}{|c|}{}\\\hline
$\gamma_{1}$ & \multicolumn{2}{|c}{$0.4$} & \multicolumn{2}{|c|}{$0.7$%
}\\\hline
\multicolumn{1}{|l|}{$n=1000$} & $\widehat{\gamma}_{1,K_{2}}$ &
$\widehat{\gamma}_{1}^{\left(  MNS\right)  }$ & $\widehat{\gamma}_{1,K_{2}}$ &
$\widehat{\gamma}_{1}^{\left(  MNS\right)  }$\\\hline
\multicolumn{5}{|c|}{$p=0.60$}\\\hline
\multicolumn{1}{|l|}{$k_{opt}$} & $234$ & $65$ & $243$ & $66$\\\hline
\multicolumn{1}{|l|}{bias} & $0.118$ & $0.122$ & $0.237$ & $0.238$\\\hline
\multicolumn{1}{|l|}{mse} & $0.016$ & $0.017$ & $0.064$ & $0.066$\\\hline
\multicolumn{5}{|c|}{$p=0.90$}\\\hline
\multicolumn{1}{|l|}{$k_{opt}$} & $249$ & $73$ & $264$ & $77$\\\hline
\multicolumn{1}{|l|}{bias} & $0.032$ & $0.053$ & $0.066$ & $0.089$\\\hline
\multicolumn{1}{|l|}{mse} & $0.002$ & $0.006$ & $0.008$ & $0.015$\\\hline
\multicolumn{5}{l|}{}%
\end{tabular}
\caption{Tail index estimation based on 2000 samples of size 1000 from Burr model right censored by another Burr model.}\label{BB}%
\end{table}%
%

\begin{table}[h] \centering
\begin{tabular}
[c]{|c|c|c|c|c|}\hline
\multicolumn{5}{|c|}{}\\\hline
$\gamma_{1}$ & \multicolumn{2}{|c}{$0.4$} & \multicolumn{2}{|c|}{$0.7$%
}\\\hline
\multicolumn{1}{|l|}{$n=1000$} & $\widehat{\gamma}_{1,K_{2}}$ &
$\widehat{\gamma}_{1}^{\left(  MNS\right)  }$ & $\widehat{\gamma}_{1,K_{2}}$ &
$\widehat{\gamma}_{1}^{\left(  MNS\right)  }$\\\hline
\multicolumn{5}{|c|}{$p=0.60$}\\\hline
\multicolumn{1}{|l|}{$k_{opt}$} & $234$ & $64$ & $234$ & $66$\\\hline
\multicolumn{1}{|l|}{bias} & $0.099$ & $0.107$ & $0.181$ & $0.189$\\\hline
\multicolumn{1}{|l|}{mse} & $0.012$ & $0.014$ & $0.039$ & $0.044$\\\hline
\multicolumn{5}{|c|}{$p=0.90$}\\\hline
\multicolumn{1}{|l|}{$k_{opt}$} & $243$ & $71$ & $243$ & $70$\\\hline
\multicolumn{1}{|l|}{bias} & $0.014$ & $0.047$ & $0.022$ & $0.080$\\\hline
\multicolumn{1}{|l|}{mse} & $0.001$ & $0.004$ & $0.004$ & $0.013$\\\hline
\multicolumn{5}{l|}{}%
\end{tabular}
\caption{Tail index estimation based on 2000 samples of size 1000 from Fr\'{e}chet model right censored  by another Fr\'{e}chet model.}\label{FF}%
\end{table}%

\section{\textbf{Real data application\label{sec4}}}

\begin{figure}[h]%
\centering
\includegraphics[
height=2.6844in,
width=5.6861in
]%
{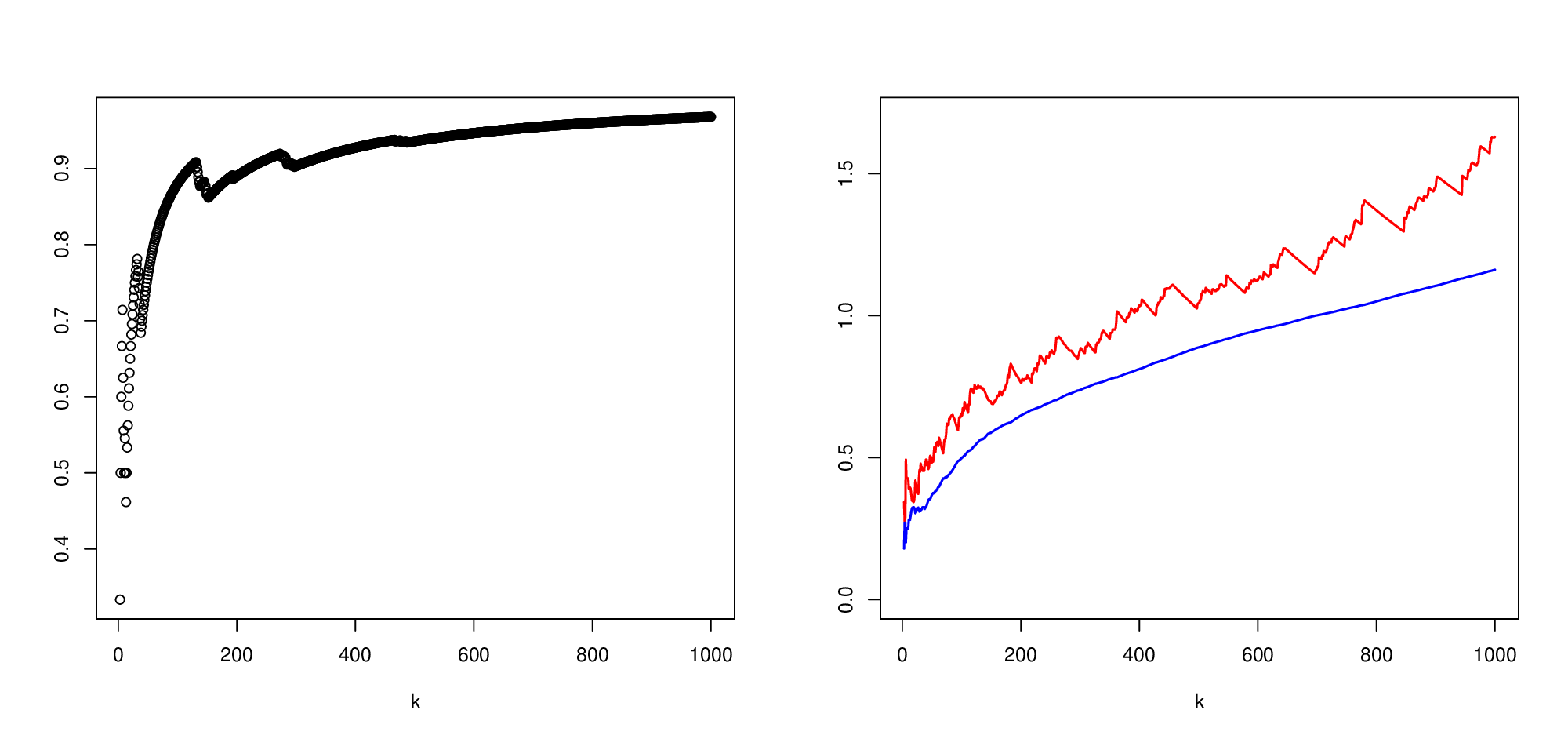}%
\caption{Insurance loss dataset. Left panel: empirical proportion of upper
non-censored observations. Right panel: $\protect\widehat{\gamma}_{1,K2}$
(blue line) and $\protect\widehat{\gamma}_{1}^{\left(  MNS\right)  }$ (red
line) as a functions of $k.$}%
\label{loss}%
\end{figure}

\noindent The data used in this study pertain to insurance losses that are
collected by the US Insurance Services Office, Inc. Available in the "copula"
package of the statistical software R, these data have been analyzed in
several studies, such as in \cite{FV98}, \cite{Klugman99} and \cite{DPV06}
amongst others. The sample consists of $1500$ observations, $34$ of which are
right censored. In insurance, censoring occurs when the loss exceeds a policy
limit representing the maximum amount a company can compensate. Each loss has
its own policy limit which varies from one contract to another. We apply our
estimation methodology to this dataset using the biweight kernel function,
resulting in the estimator $\widehat{\gamma}_{1,K_{2}}.$ These are plotted, in
the right panel of Figure \ref{loss}, as functions of the number of upper
order statistics $k$ along with $\widehat{\gamma}_{1}^{\left(  MNS\right)  }.$
We also plot the empirical proportion of upper non-censored observations
$\widehat{p}$ in the left panel of Figure \ref{loss}, where we clearly see
that $\widehat{p}>1/2$ for the vast majority of $k$ values. Our estimator
tends to provide lower values for the tail index compared to $\widehat{\gamma
}_{1}^{\left(  MNS\right)  },$ which clearly appears more volatile, especially
for small values of $k.$ Overall, the two curves exhibit an increasing trend,
reflecting the impact of larger numbers of upper order statistics on extreme
value based estimation. These results suggest that kernel estimator offer a
more smooth and stable alternative, as discussed in Section \ref{sec3}.%

\section{\textbf{Proof\label{sec5}}}

\subsection{Proof of the theorem}

\noindent We begin by proving the consistency of $\widehat{\gamma}_{1,K}.$ By
using, in formula $\left(  \ref{f-form}\right)  ,$ the change of variables
$x:=y/Z_{n-k:n}$ then by making an integration by parts, we get%
\[
\widehat{\gamma}_{1,K}=\int_{1}^{\infty}x^{-1}g_{K}\left(  \frac{\overline
{F}_{n}^{\left(  NA\right)  }\left(  Z_{n-k:n}x\right)  }{\overline{F}%
_{n}^{\left(  NA\right)  }\left(  Z_{n-k:n}\right)  }\right)  dx,
\]
where $g_{K}\left(  s\right)  :=sK\left(  s\right)  .$ Taylor's expansion to
the second order, yields that%
\begin{align*}
&  g_{K}\left(  \frac{\overline{F}_{n}^{\left(  NA\right)  }\left(
Z_{n-k:n}x\right)  }{\overline{F}_{n}^{\left(  NA\right)  }\left(
Z_{n-k:n}\right)  }\right)  -g_{K}\left(  \frac{\overline{F}\left(
Z_{n-k:n}x\right)  }{\overline{F}\left(  Z_{n-k:n}\right)  }\right) \\
&  =w_{K}\left(  \frac{\overline{F}\left(  Z_{n-k:n}x\right)  }{\overline
{F}\left(  Z_{n-k:n}\right)  }\right)  L_{n,k}+\frac{1}{2}w_{K}^{\prime
}\left(  \eta_{n}\left(  x\right)  \right)  L_{n,k}^{2},
\end{align*}
where%
\[
L_{n,k}\left(  x\right)  :=\frac{\overline{F}_{n}^{\left(  NA\right)  }\left(
Z_{n-k:n}x\right)  }{\overline{F}_{n}^{\left(  NA\right)  }\left(
Z_{n-k:n}\right)  }-\frac{\overline{F}\left(  Z_{n-k:n}x\right)  }%
{\overline{F}\left(  Z_{n-k:n}\right)  },
\]
and $\eta_{n}\left(  x\right)  $ is between $\frac{\overline{F}_{n}^{\left(
NA\right)  }\left(  Z_{n-k:n}x\right)  }{\overline{F}_{n}^{\left(  NA\right)
}\left(  Z_{n-k:n}\right)  }$ and $\frac{\overline{F}\left(  Z_{n-k:n}%
x\right)  }{\overline{F}\left(  Z_{n-k:n}\right)  }.$ The notation
$w_{K}^{\prime}\left(  s\right)  :=dw_{K}\left(  s\right)  /ds$ stands for the
Lebesgue first derivative of $w_{K}.$ Then $\widehat{\gamma}_{1,K}$ may be
splitted into the sum of
\[
I_{1n}:=\int_{1}^{\infty}x^{-1}g_{K}\left(  \frac{\overline{F}\left(
Z_{n-k:n}x\right)  }{\overline{F}\left(  Z_{n-k:n}\right)  }\right)  dx,\text{
}%
\]%
\[
I_{2n}:=\int_{1}^{\infty}x^{-1}L_{n,k}\left(  x\right)  w_{K}\left(
\frac{\overline{F}\left(  Z_{n-k:n}x\right)  }{\overline{F}\left(
Z_{n-k:n}\right)  }\right)  dx
\]
and
\[
R_{n}:=\frac{1}{2}\int_{1}^{\infty}x^{-1}L_{n,k}^{2}\left(  x\right)
w_{K}^{\prime}\left(  \eta_{n}\left(  x\right)  \right)  dx.
\]
By using an integration by parts, we show that assertion $\left(
\ref{gamma1}\right)  $ is equivalent to
\[
\int_{1}^{\infty}x^{-1}g_{K}\left(  \frac{\overline{F}\left(  ux\right)
}{\overline{F}\left(  u\right)  }\right)  dx\rightarrow\gamma_{1},\text{ as
}u\rightarrow\infty.
\]
Since $Z_{n-k:n}\overset{\mathbf{P}}{\rightarrow}\infty,$ then $I_{1n}%
\overset{\mathbf{P}}{\rightarrow}\gamma_{1}$ as $n\rightarrow\infty,$ as well.
Now, we need to show that both terms $I_{2n}$ and $R_{n}$ tend to zero in
probability. \cite{MNS2025} stated in their Gaussian approximation $\left(
6.29\right)  ,$ that there exists a sequence of standard Wiener processes
$\left\{  W_{n}\left(  s\right)  ;\text{ }0\leq s\leq1\right\}  $ defined on
the probability space $\left(  \Omega,\mathcal{A},\mathbf{P}\right)  ,$ such
that for every small $\epsilon>0$ and for any $1/4<\eta<1/2,$ we have%
\begin{equation}
\sqrt{k}L_{n,k}\left(  x\right)  =J_{n}\left(  x\right)  +o_{\mathbf{P}%
}\left(  x^{\left(  2\eta-p\right)  /\gamma+\epsilon}\right)  , \label{gauss}%
\end{equation}
uniformly over $x\geq1,$ where $J_{n}\left(  x\right)  =J_{1n}\left(
x\right)  +J_{2n}\left(  x\right)  ,$ with%
\begin{equation}
J_{1n}\left(  x\right)  :=\sqrt{\frac{n}{k}}\left\{  x^{1/\gamma_{2}%
}\mathbf{W}_{n,1}\left(  \frac{k}{n}x^{-1/\gamma}\right)  -x^{-1/\gamma_{1}%
}\mathbf{W}_{n,1}\left(  \frac{k}{n}\right)  \right\}  , \label{Jn1}%
\end{equation}
and
\begin{equation}
J_{2n}\left(  x\right)  :=\dfrac{x^{-1/\gamma_{1}}}{\gamma}\sqrt{\dfrac{n}{k}}%
{\displaystyle\int_{1}^{x}}
u^{1/\gamma-1}\left\{  p\mathbf{W}_{n,2}\left(  \dfrac{k}{n}u^{-1/\gamma
}\right)  -q\mathbf{W}_{n,1}\left(  \dfrac{k}{n}u^{-1/\gamma}\right)
\right\}  du, \label{Jn2}%
\end{equation}
where $\mathbf{W}_{n,1}$ and $\mathbf{W}_{n,2}$ are two independent Wiener
processes defined, for $0\leq s\leq1,$ by%
\[
\mathbf{W}_{n,1}\left(  s\right)  :=\left\{  W_{n}\left(  \theta\right)
-W_{n}\left(  \theta-ps\right)  \right\}  \mathbb{I}_{\left\{  \theta
-ps\geq0\right\}  }\text{ and }\mathbf{W}_{n,2}\left(  s\right)
:=W_{n}\left(  1\right)  -W_{n}\left(  1-qs\right)  ,
\]
with $\theta:=H^{\left(  1\right)  }\left(  \infty\right)  $ and
$p=1-q:=\gamma/\gamma_{1},$ uniformly on $x\geq1.$ On the other hand, since
$w_{K}$ is bounded on $\mathbb{R},$ then making use of Gaussian approximation
$\left(  \ref{gauss}\right)  ,$ we infer that
\[
\sqrt{k}\left\vert I_{2n}\right\vert =O_{\mathbf{P}}\left(  1\right)  \int%
_{1}^{\infty}x^{-1}\left\vert J_{n}\left(  x\right)  \right\vert
dx+o_{\mathbf{P}}\left(  1\right)  \int_{1}^{\infty}x^{\left(  2\eta-p\right)
/\gamma+\epsilon-1}dx.
\]
Note that $p>1/2$ (see Remark $\left(  \ref{remark4}\right)  $), then the
integral $\int_{1}^{\infty}x^{\left(  2\eta-p\right)  /\gamma+\epsilon-1}dx$
is finite for any small $\epsilon>0$ and every $1/4<\eta<1/2.$ Next we give an
approximation to $J_{n}\left(  x\right)  .$ From Lemma 3.2 in \cite{EHL2006}
\[
\sup_{0<s\leq1}\frac{\left\vert W_{n}\left(  s\right)  \right\vert
}{s^{1-2\eta}}=O_{\mathbf{P}}\left(  1\right)  ,\text{ }%
\]
for every $0<\eta<1/4.$ On the other hand, for each $n$
\[
\mathbf{W}_{n,1}\left(  s\right)  \overset{\mathcal{D}}{=}W_{n}\left(
ps\right)  \overset{\mathcal{D}}{=}p^{1/2}W_{n}\left(  s\right)  \text{ and
}\mathbf{W}_{n,2}\left(  s\right)  \overset{\mathcal{D}}{=}W_{n}\left(
ps\right)  \overset{\mathcal{D}}{=}q^{1/2}W_{n}\left(  s\right)  ,
\]
jointly on $0<s\leq1.$ The notation $\overset{\mathcal{D}}{=}$ stands for the
equality in distribution. It follows that
\begin{equation}
\sup_{0<s\leq1}\frac{\left\vert \mathbf{W}_{n,i}\left(  s\right)  \right\vert
}{s^{1-2\eta}}=O_{\mathbf{P}}\left(  1\right)  ,\text{ }i=1,2. \label{sup}%
\end{equation}
It is clear that, for each $n$%
\[
J_{1n}\left(  x\right)  \overset{\mathcal{D}}{=}x^{1/\gamma_{2}}%
\mathbf{W}_{n,1}\left(  x^{-1/\gamma}\right)  -x^{-1/\gamma_{1}}%
\mathbf{W}_{n,1}\left(  1\right)  ,\text{ jointly on }x\geq1.
\]
Recall that $\gamma_{2}=\left(  1-p\right)  /\gamma$ and $\gamma_{1}%
=p/\gamma,$ then without loss of generality, we may write%
\[
J_{1n}\left(  x\right)  =x^{\left(  1-p\right)  /\gamma}\mathbf{W}%
_{n,1}\left(  x^{-1/\gamma}\right)  -x^{-p/\gamma}\mathbf{W}_{n,1}\left(
1\right)  .
\]
By applying assertion $\left(  \ref{sup}\right)  ,$ for $i=1,$ it is easy to
verify that $J_{1n}\left(  x\right)  =O_{\mathbf{P}}\left(  x^{\left(
2\eta-p\right)  /\gamma}\right)  ,$ uniformly over $x\geq1.$ Likewise, without
loss of generality, we may also write%
\[
J_{2n}\left(  x\right)  =\dfrac{x^{-p/\gamma}}{\gamma}%
{\displaystyle\int_{1}^{x}}
u^{1/\gamma-1}\left\{  p\mathbf{W}_{n,2}\left(  u^{-1/\gamma}\right)
+q\mathbf{W}_{n,1}\left(  u^{-1/\gamma}\right)  \right\}  du.
\]
On again by applying assertion $\left(  \ref{sup}\right)  ,$ for $i=2,$ we
also show that $J_{1n}\left(  x\right)  =O_{\mathbf{P}}\left(  x^{\left(
2\eta-p\right)  /\gamma}\right)  ,$ uniformly over $x\geq1.$ In summary, we
showed that for ever $1/4<\eta<1/2:$
\begin{equation}
J_{n}\left(  x\right)  =O_{\mathbf{P}}\left(  x^{\left(  2\eta-p\right)
/\gamma}\right)  ,\text{ uniformly over }x\geq1. \label{Jn}%
\end{equation}
It is abvious that $\left(  \ref{Jn}\right)  $ implies that
\[
\int_{1}^{\infty}x^{-1}\left\vert J_{n}\left(  x\right)  \right\vert
dx=O_{\mathbf{P}}\left(  1\right)  \int_{1}^{\infty}x^{-1+\left(
2\eta-p\right)  /\gamma}dx=O_{\mathbf{P}}\left(  1\right)  ,
\]
provided that $p>2\eta<1/2.$ By considering the latter assumption, we infer
that $I_{2n}\overset{\mathbf{P}}{\rightarrow}0$ as $n\rightarrow\infty$
(because $1/\sqrt{k}\rightarrow0).$ Let us now consider the remainder term
$R_{n}.$ Let us write%
\[
kR_{n}=\frac{1}{2}\int_{1}^{\infty}x^{-1}\left(  \sqrt{k}L_{n,k}\left(
x\right)  \right)  ^{2}w_{K}^{\prime}\left(  \eta_{n}\left(  x\right)
\right)  dx.
\]
Once again, by using $\left(  \ref{gauss}\right)  $ with the fact that
$w_{K}^{\prime}$ is bounded, we get
\begin{align*}
kR_{n}  &  =O_{\mathbf{P}}\left(  1\right)  \int_{1}^{\infty}x^{-1}J_{n}%
^{2}\left(  x\right)  dx\\
&  +o_{\mathbf{P}}\left(  1\right)  \int_{1}^{\infty}x^{\left(  2\eta
-p\right)  /\gamma+\epsilon-1}\left\vert J_{n}\left(  x\right)  \right\vert
dx+\int_{1}^{\infty}o_{\mathbf{P}}\left(  x^{2\left(  2\eta-p\right)
/\gamma+2\epsilon-1}\right)  dx.
\end{align*}
It is clear that both the second and the third terms in the right-hand side
above equal to $o_{\mathbf{P}}\left(  1\right)  .$ Once agrain, by considering
assumption $p>1/2$ and approximation $\left(  \ref{Jn}\right)  ,$ we write
\[
\int_{1}^{\infty}x^{-1}J_{n}^{2}\left(  x\right)  dx=O_{\mathbf{P}}\left(
1\right)  \int_{1}^{\infty}x^{-1+2\left(  2\eta-p\right)  /\gamma
}dx=O_{\mathbf{P}}\left(  1\right)  ,
\]
therefore $R_{n}\overset{\mathbf{P}}{\rightarrow}0$ when $n\rightarrow\infty.$
In summary, we showed that both terms $I_{2n}$ and $R_{n}$ tend in probability
to zero, while $I_{1n}\overset{\mathbf{P}}{\rightarrow}\gamma_{1},$ this means
that $\widehat{\gamma}_{1,K}\overset{\mathbf{P}}{\rightarrow}\gamma_{1}$ as
$n\rightarrow\infty.$ Let us now establish the asymptotic normality of
$\widehat{\gamma}_{1,K}.$ Recall that $\int_{\mathbb{R}}K\left(  t\right)
dt=1,$ then it is easy to verify that $\int_{1}^{\infty}x^{-1}g_{K}\left(
x^{-1/\gamma_{1}}\right)  dx=\gamma_{1},$ which allows us to write%
\[
\widehat{\gamma}_{1,K}-\gamma_{1}=\int_{1}^{\infty}x^{-1}\left\{  g_{K}\left(
\frac{\overline{F}_{n}^{\left(  NA\right)  }\left(  Z_{n-k:n}x\right)
}{\overline{F}_{n}^{\left(  NA\right)  }\left(  Z_{n-k:n}\right)  }\right)
-g_{K}\left(  x^{-1/\gamma_{1}}\right)  \right\}  dx.
\]
Once again, by using the second-order Taylor's expansion to $g_{K},$\ yields%
\begin{align*}
&  \sqrt{k}\left(  \widehat{\gamma}_{1,K}-\gamma_{1}\right) \\
&  =\int_{1}^{\infty}x^{-1}D_{n}\left(  x\right)  w_{K}\left(  x^{-1/\gamma
_{1}}\right)  dx+\frac{1}{2\sqrt{k}}\int_{1}^{\infty}x^{-1}D_{n}^{2}\left(
x\right)  w_{K}^{\prime}\left(  \zeta_{n}\left(  x\right)  \right)  dx,
\end{align*}
where%
\[
D_{n}\left(  x\right)  :=\sqrt{k}\left(  \frac{\overline{F}_{n}^{\left(
NA\right)  }\left(  Z_{n-k:n}x\right)  }{\overline{F}_{n}^{\left(  NA\right)
}\left(  Z_{n-k:n}\right)  }-x^{-1/\gamma_{1}}\right)  ,
\]
with $\zeta_{n}\left(  x\right)  $\ being between $\overline{F}_{n}^{\left(
NA\right)  }\left(  xZ_{n-k:n}\right)  /\overline{F}_{n}^{\left(  NA\right)
}\left(  Z_{n-k:n}\right)  $ and $x^{1/\gamma_{1}}.$ \cite{MNS2025} stated in
their Theorem 3, that%
\begin{equation}
\sup_{x\geq p^{\gamma}}x^{\epsilon/p\gamma_{1}}\left\vert D_{n}\left(
x\right)  -J_{n}\left(  x\right)  -x^{-1/\gamma_{1}}\dfrac{x^{\tau_{1}%
/\gamma_{1}}-1}{\tau_{1}\gamma_{1}}\sqrt{k}A_{1}\left(  h\right)  \right\vert
\overset{\mathbf{P}}{\rightarrow}0, \label{approx}%
\end{equation}
for every $0<\epsilon<1/2,$ as $n\rightarrow\infty,$ provided that $p>1/2$ and
$\sqrt{k}A_{1}\left(  h\right)  =O\left(  1\right)  .$ In order to establish
the asymptotic normality of $\widehat{\gamma}_{1,K},$ we will use the
following decomposition%
\[
\sqrt{k}\left(  \widehat{\gamma}_{1,K}-\gamma_{1}\right)  =\int_{1}^{\infty
}x^{-1}J_{n}\left(  x\right)  w_{K}\left(  x^{-1/\gamma_{1}}\right)
dx+\sum_{i=1}^{3}\mathcal{R}_{ni},
\]
where%
\[
\mathcal{R}_{n1}:=\frac{1}{2\sqrt{k}}\int_{1}^{\infty}x^{-1}D_{n}^{2}\left(
x\right)  w_{K}^{\prime}\left(  \zeta_{n}\left(  x\right)  \right)  dx,
\]%
\[
\mathcal{R}_{n2}:=\int_{1}^{\infty}x^{-1}\left\{  D_{n}\left(  x\right)
-J_{n}\left(  x\right)  -x^{-1/\gamma_{1}}\dfrac{x^{\tau_{1}/\gamma_{1}}%
-1}{\gamma_{1}\tau_{1}}\sqrt{k}A_{1}\left(  h\right)  \right\}  w_{K}\left(
x^{-1/\gamma_{1}}\right)  dx,
\]
and
\[
\mathcal{R}_{n3}:=\sqrt{k}A_{1}\left(  h\right)  \int_{1}^{\infty}%
x^{-1}\left\{  x^{-1/\gamma_{1}}\dfrac{x^{\tau_{1}/\gamma_{1}}-1}{\gamma
_{1}\tau_{1}}\right\}  w_{K}\left(  x^{-1/\gamma_{1}}\right)  dx.
\]
Next, we show that $\mathcal{R}_{ni}=o_{\mathbf{P}}\left(  1\right)  ,$ for
$i=1,2$ as $n\rightarrow\infty.$ Indeed, since $w_{K}^{\prime}$ is bounded
then
\[
\int_{1}^{\infty}x^{-1}D_{n}^{2}\left(  x\right)  w_{K}^{\prime}\left(
\zeta_{n}\left(  x\right)  \right)  dx=O_{\mathbf{P}}\left(  1\right)
\int_{1}^{\infty}x^{-1}D_{n}^{2}\left(  x\right)  dx.
\]
Next we show that $\int_{1}^{\infty}x^{-1}D_{n}^{2}\left(  x\right)
dx=O_{\mathbf{P}}\left(  1\right)  .$ In view of weak approximation $\left(
\ref{approx}\right)  ,$ we write
\[
\int_{1}^{\infty}x^{-1}D_{n}^{2}\left(  x\right)  dx=\int_{1}^{\infty}%
x^{-1}\left(  J_{n}\left(  x\right)  -x^{-1/\gamma_{1}}\dfrac{x^{\tau
_{1}/\gamma_{1}}-1}{\tau_{1}\gamma_{1}}\sqrt{k}A_{1}\left(  h\right)
+o_{\mathbf{P}}\left(  x^{-\epsilon/\left(  p\gamma_{1}\right)  }\right)
\right)  ^{2}dx.
\]
Once again, making use of approximation $\left(  \ref{sup}\right)  $ and using
the fact that $\sqrt{k}A_{1}\left(  h\right)  =O\left(  1\right)  ,$ we show
that the right-hand side of the equation above is asymptotically bounded in
probability, therefore we omit the details. It follows that $\mathcal{R}%
_{ni}=O_{\mathbf{P}}\left(  k^{-1/2}\right)  $ which tends to zero in
probability as $k^{-1/2}\rightarrow0.$ Likewise since $w_{k}$ is bounded, then
(once again) by using $\left(  \ref{approx}\right)  ,$ we get $\mathcal{R}%
_{n2}=o_{\mathbf{P}}\left(  1/\sqrt{k}\right)  .$ The change of variables
$s=x^{-1/\gamma_{1}}$ and an integration by parts transform the previous
integral into%
\[
\int_{1}^{\infty}x^{-1}\left\{  x^{-1/\gamma_{1}}\dfrac{x^{\tau_{1}/\gamma
_{1}}-1}{\gamma_{1}\tau_{1}}\right\}  w_{K}\left(  x^{-1/\gamma_{1}}\right)
dx=\int_{0}^{1}s^{-\tau_{1}}K\left(  s\right)  ds.
\]
It follows that
\[
\mathcal{R}_{n3}=\sqrt{k}A_{1}\left(  h\right)  \int_{0}^{1}s^{-\tau_{1}%
}K\left(  s\right)  ds,\text{ as }n\rightarrow\infty.
\]
In summary, we demonstrated that
\[
\sqrt{k}\left(  \widehat{\gamma}_{1,K}-\gamma_{1}\right)  =T_{n}+\sqrt{k}%
A_{1}\left(  h\right)  \int_{0}^{1}s^{-\tau_{1}}K\left(  s\right)
ds+o_{\mathbf{P}}\left(  1\right)  ,\text{ as }n\rightarrow\infty,
\]
where $T_{n}:=\int_{1}^{\infty}x^{-1}J_{n}\left(  x\right)  w_{K}\left(
x^{-1/\gamma_{1}}\right)  dx$ is a sequence of centered Gaussian rv's, whose
variancce will be computed right below.\medskip

\noindent\textbf{Computation of }$\sigma_{K}^{2}$\smallskip

\noindent The change of variables $t=x^{-1/\gamma_{1}}$ with the definition of
$w_{K}$ yield that%
\[
T_{n}=\gamma_{1}\int_{0}^{1}t^{-1}J_{n}\left(  t^{-\gamma_{1}}\right)
d\left\{  tK\left(  t\right)  \right\}  ,
\]
where $J_{n}\left(  t^{-\gamma_{1}}\right)  =J_{1n}\left(  t^{-\gamma_{1}%
}\right)  +J_{2n}\left(  t^{-\gamma_{1}}\right)  ,$ with $J_{1n}$ and
$J_{2n\text{ }}$ being defined in $\left(  \ref{Jn1}\right)  $ and $\left(
\ref{Jn2}\right)  $\ respectively. Note that, as a linear functional of two
(independent) sequences of Wiener processes $\mathbf{W}_{n,1}$ and
$\mathbf{W}_{n,2},$ $T_{n}$ is a sequence of centred Gaussian rv's. Now, we
write $T_{n}$ into the following sum%
\[
T_{n}=\gamma_{1}\int_{0}^{1}t^{-1}J_{n1}\left(  t^{-\gamma_{1}}\right)
d\left\{  tK\left(  t\right)  \right\}  +\gamma_{1}\int_{0}^{1}t^{-1}%
J_{n2}\left(  t^{-\gamma_{1}}\right)  d\left\{  tK\left(  t\right)  \right\}
.
\]
In other words, by squaring $T_{n},$ we may write $\mathbf{E}\left[  T_{n}%
^{2}\right]  $ as the sum of the following three terms%
\[
S_{1}:=\gamma_{1}^{2}\mathbf{E}\left(  \int_{0}^{1}t^{-1}J_{n1}\left(
t^{-\gamma_{1}}\right)  d\left\{  tK\left(  t\right)  \right\}  \right)
^{2},
\]%
\[
S_{2}:=\gamma_{1}^{2}\mathbf{E}\left(  \int_{0}^{1}t^{-1}J_{n2}\left(
t^{-\gamma_{1}}\right)  d\left\{  tK\left(  t\right)  \right\}  \right)
^{2},
\]
and%
\[
S_{3}:=2\gamma_{1}^{2}\mathbf{E}\left(  \int_{0}^{1}t^{-1}J_{n1}\left(
t^{-\gamma_{1}}\right)  d\left\{  tK\left(  t\right)  \right\}  \int_{0}%
^{1}t^{-1}J_{n2}\left(  t^{-\gamma_{1}}\right)  d\left\{  tK\left(  t\right)
\right\}  \right)  ,
\]
for the computaion of which the following covariances
\begin{equation}
\mathbf{E}\left[  \mathbf{W}_{n,1}\left(  s\right)  \mathbf{W}_{n,2}\left(
t\right)  \right]  =0, \label{W1W2}%
\end{equation}%
\begin{equation}
\mathbf{E}\left[  \mathbf{W}_{n,1}\left(  s\right)  \mathbf{W}_{n,1}\left(
t\right)  \right]  =p\min\left(  s,t\right)  \label{W1W1}%
\end{equation}
and%
\begin{equation}
\mathbf{E}\left[  \mathbf{W}_{n,2}\left(  s\right)  \mathbf{W}_{n,2}\left(
t\right)  \right]  =q\min\left(  s,t\right)  , \label{W2W2}%
\end{equation}
will be useful. We start with the calculation of $S_{1}.$ We have%
\[
S_{1}=\gamma_{1}^{2}\int_{0}^{1}\int_{0}^{1}s^{-1}t^{-1}\mathbf{E}\left[
J_{n1}\left(  s^{-\gamma_{1}}\right)  J_{n1}\left(  t^{-\gamma_{1}}\right)
\right]  d\left\{  sK\left(  s\right)  \right\}  \left\{  tK\left(  t\right)
\right\}  .
\]
From representation $\left(  \ref{Jn1}\right)  $\ and using $\left(
\ref{W1W1}\right)  ,$ it can be readily checked that%
\[
S_{1}=2\gamma_{1}^{2}p\int_{0}^{1}s^{-\gamma_{1}/\gamma_{2}-1}\left\{
sK\left(  s\right)  \right\}  d\left\{  sK\left(  s\right)  \right\}
=\gamma_{1}^{2}p\int_{0}^{1}s^{-\gamma_{1}/\gamma_{2}-1}d\left\{  sK\left(
s\right)  \right\}  ^{2}.
\]
Note that $-\gamma_{1}/\gamma_{2}-1=-\gamma_{1}/\gamma,$ then it follows that
$S_{1}=\gamma_{1}^{2}p\int_{0}^{1}s^{-\gamma_{1}/\gamma}d\left\{  sK\left(
s\right)  \right\}  ^{2}.$ To compute $S_{2}$ and $S_{3},$ we need to modify
$J_{2n}\left(  x\right)  ,$ defined in $\left(  \ref{Jn2}\right)  .$ We do
this by a change of variables and we get%
\begin{equation}
J_{2n}\left(  x\right)  =-x^{-1/\gamma_{1}}\sqrt{\dfrac{n}{k}}%
{\displaystyle\int_{1}^{x^{-1/\gamma}}}
s^{-2}\left\{  p\mathbf{W}_{n,2}\left(  \dfrac{k}{n}s\right)  -q\mathbf{W}%
_{n,1}\left(  \dfrac{k}{n}s\right)  \right\}  ds. \label{J2n}%
\end{equation}
Note that $S_{2}=\mathbf{E}\left[  I^{2}\right]  ,$ with $I:=\gamma_{1}%
\int_{0}^{1}t^{-1}J_{n2}\left(  t^{-\gamma_{1}}\right)  d\left\{  tK\left(
t\right)  \right\}  ,$ which we simplify by using $\left(  \ref{J2n}\right)  $
and after an integration by parts, we get
\[
I=\frac{\gamma_{1}^{2}}{\gamma}\sqrt{\dfrac{n}{k}}\int_{0}^{1}K\left(
t\right)  t^{-\gamma_{1}/\gamma}\left[  p\mathbf{W}_{n,2}\left(  \dfrac{k}%
{n}t^{\gamma_{1}/\gamma}\right)  -q\mathbf{W}_{n,1}\left(  \dfrac{k}%
{n}t^{\gamma_{1}/\gamma}\right)  \right]  dt.
\]
By using results $\left(  \ref{W1W2}\right)  ,$ $\left(  \ref{W1W1}\right)  $
and $\left(  \ref{W2W2}\right)  ,$\ we obtain%
\[
S_{2}=\left(  \gamma_{1}^{3}/\gamma\right)  q\int_{0}^{1}K\left(  s\right)
s^{-\gamma_{1}/\gamma}\int_{0}^{s}K\left(  t\right)  dtds+\int_{0}^{1}\int%
_{s}^{1}K\left(  t\right)  t^{-\gamma_{1}/\gamma}dtd\left\{  \int_{0}%
^{s}K\left(  t\right)  dt\right\}  .
\]
Since $\gamma_{1}/\gamma=1/p>1/2$ and after an integration by parts, we infer
that%
\[
S_{2}=2\gamma_{1}^{2}\frac{q}{p}\int_{0}^{1}\left\{  \int_{0}^{s}K\left(
t\right)  dt\right\}  K\left(  s\right)  s^{-\gamma_{1}/\gamma}ds.
\]
Finally, we turn our attention to $S_{3}=2\gamma_{1}^{2}\mathbf{E}\left[
I^{\prime}\right]  ,$ where%
\[
I^{\prime}:=\int_{0}^{1}t^{-1}J_{n1}\left(  t^{-\gamma_{1}}\right)  d\left\{
tK\left(  t\right)  \right\}  \int_{0}^{1}t^{-1}J_{n2}\left(  t^{-\gamma_{1}%
}\right)  d\left\{  tK\left(  t\right)  \right\}  .
\]
By substituting equations $\left(  \ref{Jn1}\right)  $ and $\left(
\ref{J2n}\right)  $ in $I^{\prime},$ we have%
\begin{align*}
I^{\prime}  &  =-\frac{n}{k}\int_{0}^{1}\left\{  t^{-\gamma_{1}/\gamma_{2}%
-1}\mathbf{W}_{n,1}\left(  \frac{k}{n}t^{\gamma_{1}/\gamma}\right)
-\mathbf{W}_{n,1}\left(  \frac{k}{n}\right)  \right\}  d\left\{  tK\left(
t\right)  \right\} \\
&  \times\int_{0}^{1}%
{\displaystyle\int_{1}^{t^{\gamma_{1}/\gamma}}}
s^{-2}\left\{  p\mathbf{W}_{n,2}\left(  \dfrac{k}{n}s\right)  -q\mathbf{W}%
_{n,1}\left(  \dfrac{k}{n}s\right)  \right\}  dsd\left\{  tK\left(  t\right)
\right\}  .
\end{align*}
Integrating the second factor in $I^{\prime}$ by parts yields%
\begin{align*}
I^{\prime}  &  =\frac{n}{k}\frac{\gamma_{1}}{\gamma}\int_{0}^{1}\left\{
t^{-\gamma_{1}/\gamma_{2}-1}\mathbf{W}_{n,1}\left(  \frac{k}{n}t^{\gamma
_{1}/\gamma}\right)  -\mathbf{W}_{n,1}\left(  \frac{k}{n}\right)  \right\}
d\left\{  tK\left(  t\right)  \right\} \\
&  \times\int_{0}^{1}K\left(  t\right)  t^{-\gamma_{1}/\gamma}\left\{
p\mathbf{W}_{n,2}\left(  \dfrac{k}{n}t^{\gamma_{1}/\gamma}\right)
-q\mathbf{W}_{n,1}\left(  \dfrac{k}{n}t^{\gamma_{1}/\gamma}\right)  \right\}
dt,
\end{align*}
which may be decomposed into the sum of
\[
I_{1}^{\prime}:=\frac{n}{k}\int_{0}^{1}\int_{0}^{1}s^{-\gamma_{1}/\gamma
_{2}-1}\mathbf{W}_{n,1}\left(  \frac{k}{n}s^{\gamma_{1}/\gamma}\right)
\mathbf{W}_{n,2}\left(  \dfrac{k}{n}t^{\gamma_{1}/\gamma}\right)  K\left(
t\right)  t^{-\gamma_{1}/\gamma}d\left\{  sK\left(  s\right)  \right\}  dt,
\]%
\[
I_{2}^{\prime}:=-q\frac{\gamma_{1}}{\gamma}\frac{n}{k}\int_{0}^{1}\int_{0}%
^{1}s^{-\gamma_{1}/\gamma_{2}-1}\mathbf{W}_{n,1}\left(  \frac{k}{n}%
s^{\gamma_{1}/\gamma}\right)  \mathbf{W}_{n,1}\left(  \dfrac{k}{n}%
t^{\gamma_{1}/\gamma}\right)  K\left(  t\right)  t^{-\gamma_{1}/\gamma
}d\left\{  sK\left(  s\right)  \right\}  dt,
\]%
\[
I_{3}^{\prime}:=-\frac{n}{k}\int_{0}^{1}\int_{0}^{1}\mathbf{W}_{n,1}\left(
\frac{k}{n}\right)  \mathbf{W}_{n,2}\left(  \dfrac{k}{n}t^{\gamma_{1}/\gamma
}\right)  K\left(  t\right)  t^{-\gamma_{1}/\gamma}d\left\{  sK\left(
s\right)  \right\}  dt,
\]
and
\[
I_{4}^{\prime}:=q\frac{\gamma_{1}}{\gamma}\frac{n}{k}\int_{0}^{1}\int_{0}%
^{1}\mathbf{W}_{n,1}\left(  \frac{k}{n}\right)  \mathbf{W}_{n,1}\left(
\dfrac{k}{n}t^{\gamma_{1}/\gamma}\right)  K\left(  t\right)  t^{-\gamma
_{1}/\gamma}d\left\{  sK\left(  s\right)  \right\}  dt.
\]
It is clear that, from $\left(  \ref{W1W2}\right)  ,$\ we have $\mathbf{E}%
\left[  I_{1}^{\prime}\right]  =\mathbf{E}\left[  I_{3}^{\prime}\right]  =0$
and by using $\left(  \ref{W1W1}\right)  ,$ we can readily show that
$\mathbf{E}\left[  I_{4}^{\prime}\right]  $ is zero as well. Hence, we have
$\mathbf{E}\left[  I^{\prime}\right]  =\mathbf{E}\left[  I_{2}^{\prime
}\right]  ,$ to which we apply $\left(  \ref{W1W1}\right)  $ to get%
\begin{align*}
\mathbf{E}\left[  I_{2}^{\prime}\right]   &  =-q\int_{0}^{1}\int_{0}%
^{1}s^{-\gamma_{1}/\gamma_{2}-1}\min\left(  s^{\gamma_{1}/\gamma}%
,t^{\gamma_{1}/\gamma}\right)  t^{-\gamma_{1}/\gamma}K\left(  t\right)
d\left\{  sK\left(  s\right)  \right\}  dt\\
&  =-q\int_{0}^{1}\left\{  \int_{0}^{s}K\left(  t\right)  dt\right\}
s^{-\gamma_{1}/\gamma_{2}-1}d\left\{  sK\left(  s\right)  \right\} \\
&  -q\int_{0}^{1}\left\{  \int_{s}^{1}t^{-\gamma_{1}/\gamma}K\left(  t\right)
dt\right\}  s^{-\gamma_{1}/\gamma_{2}+\gamma_{1}/\gamma-1}d\left\{  sK\left(
s\right)  \right\}  .
\end{align*}
Note that $-\gamma_{1}/\gamma_{2}+\gamma_{1}/\gamma-1=0,$ then $\mathbf{E}%
\left[  I_{2}^{\prime}\right]  =-q\left(  I_{21}^{\prime}+I_{22}^{\prime
}\right)  ,$ where
\[
I_{21}^{\prime}:=\int_{0}^{1}\left\{  \int_{0}^{s}K\left(  t\right)
dt\right\}  s^{-\gamma_{1}/\gamma_{2}-1}d\left\{  sK\left(  s\right)
\right\}  ,
\]
and
\[
I_{22}^{\prime}:=\int_{0}^{1}\left\{  \int_{s}^{1}t^{-\gamma_{1}/\gamma
}K\left(  t\right)  dt\right\}  d\left\{  sK\left(  s\right)  \right\}  ,
\]
which after integrations by parts become%
\[
I_{21}^{\prime}=-\int_{0}^{1}\left\{  K\left(  s\right)  \right\}
^{2}s^{-\gamma_{1}/\gamma_{2}}ds+\frac{\gamma_{1}}{\gamma}\int_{0}^{1}K\left(
s\right)  s^{-\gamma_{1}/\gamma_{2}-1}\left\{  \int_{0}^{s}K\left(  t\right)
dt\right\}  ds,
\]
and $I_{22}^{\prime}=\int_{0}^{1}\left(  K\left(  s\right)  \right)
^{2}s^{-\gamma_{1}/\gamma_{2}}ds.$ Thus we, obtain%
\begin{align*}
\mathbf{E}\left[  I_{2}^{\prime}\right]   &  =q\int_{0}^{1}\left(  K\left(
s\right)  \right)  ^{2}s^{-\gamma_{1}/\gamma_{2}}ds-q\int_{0}^{1}\left(
K\left(  s\right)  \right)  ^{2}s^{-\gamma_{1}/\gamma_{2}}ds\\
&  -q\frac{\gamma_{1}}{\gamma}\int_{0}^{1}\left\{  \int_{0}^{s}K\left(
t\right)  dt\right\}  K\left(  s\right)  s^{-\gamma_{1}/\gamma_{2}-1}ds.
\end{align*}
Finally, we end up with
\[
\mathbf{E}\left[  T_{n}^{2}\right]  =2\gamma_{1}^{2}p\int_{0}^{1}%
s^{-\gamma_{1}/\gamma}\left\{  sK\left(  s\right)  \right\}  d\left\{
sK\left(  s\right)  \right\}  =\gamma_{1}^{2}p\int_{0}^{1}s^{-1/p}d\left\{
sK\left(  s\right)  \right\}  ^{2},
\]
which, by an integration by parts, becomes $\gamma_{1}^{2}\int_{0}%
^{1}s^{-1/p+1}K^{2}\left(  s\right)  ds$ which is the expression of
$\sigma_{K}^{2}$ given in the theorem. Consequently, the sequence $T_{n}$
equals in distribution to $\mathcal{N}\left(  0,\sigma_{K}^{2}\right)  $ and
therefore we have%
\[
\sqrt{k}\left(  \widehat{\gamma}_{1,K}-\gamma_{1}\right)  \overset{\mathcal{D}%
}{=}\mathcal{N}\left(  0,\sigma_{K}^{2}\right)  +\sqrt{k}A_{1}\left(
h\right)  \mu_{K}+o_{\mathbf{P}}\left(  1\right)  ,
\]
which meets approximation $\left(  \ref{approx1}\right)  .$ If, in addition,
we assune that $\sqrt{k}A_{1}\left(  h\right)  \rightarrow\lambda<\infty,$ we
get%
\[
\sqrt{k}\left(  \widehat{\gamma}_{1,K}-\gamma_{1}\right)  \overset{\mathcal{D}%
}{\rightarrow}\mathcal{N}\left(  \lambda\mu_{K},\sigma_{K}^{2}\right)  ,\text{
as }n\rightarrow\infty.
\]
This completes the proof of the theorem.

\subsection{Proof of the corollary}

We saw, in Section \ref{kopt}, that the optimal sampe fraction $k_{opt}$ is
the $k$-value that minimizes the AMSE of $\widehat{\gamma}_{1,K}.$ That is%
\[
k_{opt}:=\arg\min_{1<k<n}AMSE\left(  k\right)  ,
\]
where $AMSE\left(  k\right)  $ is defined later on in $\left(  \ref{mse}%
\right)  .$ From statement $\left(  \ref{approx1}\right)  ,$ we have
\begin{align*}
&  \widehat{\gamma}_{1,K}-\gamma_{1}\\
&  \overset{\mathcal{D}}{=}\frac{\gamma_{1}}{\sqrt{k}}\mathcal{N}\left(
0,\int_{0}^{1}s^{-1/p+1}K^{2}\left(  s\right)  ds\right)  +\left(
1+o_{\mathbf{P}}\left(  1\right)  \right)  A_{1}\left(  h\right)  \int_{0}%
^{1}s^{-\tau_{1}}K\left(  s\right)  ds,
\end{align*}
as $n\rightarrow\infty,$ where $A_{1}\left(  h\right)  =A_{1}^{\ast}\left(
1/\overline{F}\left(  h\right)  \right)  ,$ with $h=U_{H}\left(  n/k\right)
.$ The corresponding AMSE of $\widehat{\gamma}_{1,K}$ equals%
\begin{equation}
AMSE\left(  k\right)  :=\left\{  \frac{\gamma_{1}^{2}}{k}\int_{0}%
^{1}s^{-1/p+1}K^{2}\left(  s\right)  ds\right\}  +\left\{  A_{1}^{2}\left(
h\right)  \left(  \int_{0}^{1}s^{-\tau_{1}}K\left(  s\right)  ds\right)
^{2}\right\}  . \label{mse}%
\end{equation}
To determine $k_{opt},$ we use similar arguments as those made with the
classical Hill estimator, see \cite{deHF06} in page 76, where the authors
consider the following form to Hall's model (statement (3.2.9)):%
\[
\overline{F}\left(  x\right)  =c_{1}x^{-1/\gamma_{1}}+c_{2}x^{-1/\gamma
_{1}+\tau_{1}/\gamma_{1}}+o\left(  x^{-1/\gamma_{1}+\tau_{1}/\gamma_{1}%
}\right)  ,\text{ as }x\rightarrow\infty.
\]
Identifying the above with $\left(  \ref{HallF}\right)  $ yields $c_{1}%
:=C_{1}$ and $c_{2}:=C_{1}C_{2}.$ Let us now give the asymptotic expression of
$A_{1}\left(  h\right)  .$ By using Hall's model $\left(  \ref{HallG}\right)
,$ we obtain
\[
\overline{H}\left(  x\right)  =\overline{F}\left(  x\right)  \overline
{G}\left(  x\right)  =\left(  1+o\left(  1\right)  \right)  C_{1}%
D_{1}x^{-1/\gamma},\text{ as }x\rightarrow\infty,
\]
where $1/\gamma=1/\gamma_{1}+1/\gamma_{2},$ which is equivalent to%
\begin{equation}
U_{H}\left(  t\right)  =\left(  1+o\left(  1\right)  \right)  \left(
C_{1}D_{1}\right)  ^{\gamma}t^{\gamma},\text{ as }t\rightarrow\infty.
\label{UH}%
\end{equation}
Combining statements $\left(  \ref{A1}\right)  $ and $\left(  \ref{UH}\right)
$ yields
\[
A_{1}\left(  h\right)  =A_{\ast}\left(  1/\overline{F}\left(  U_{H}\left(
n/k\right)  \right)  \right)  =\frac{\tau_{1}}{\gamma_{1}}C_{2}\left(  \left(
C_{1}D_{1}\right)  ^{\gamma}\left(  n/k\right)  ^{\gamma}\right)  ^{\tau
_{1}/\gamma_{1}}=c\left(  n/k\right)  ^{p\tau_{1}},
\]
where $c:=\frac{\tau_{1}}{\gamma_{1}}C_{2}\left(  C_{1}D_{1}\right)
^{p\tau_{1}}.$ Thus expression $\left(  \ref{mse}\right)  $ becomes
\[
AMSE\left(  k\right)  =\frac{\gamma_{1}^{2}}{k}\int_{0}^{1}s^{-1/p+1}%
K^{2}\left(  s\right)  ds+c^{2}\left(  \frac{n}{k}\right)  ^{2p\tau_{1}%
}\left(  \int_{0}^{1}s^{-\tau_{1}}K\left(  s\right)  ds\right)  ^{2}.
\]
By similar arguments as those used in \cite{deHF06}, page 78, we end up with
$k_{opt}=\left[  \omega n^{-2p\tau_{1}/\left(  1-2p\tau_{1}\right)  }\right]
,$ where%
\[
\omega:=\left\{  \frac{\gamma_{1}^{2}}{-2\tau_{1}c^{2}}\left(  \int_{0}%
^{1}s^{-1/p+1}K^{2}\left(  s\right)  ds\right)  \left(  \int_{0}^{1}%
s^{-\tau_{1}}K\left(  s\right)  ds\right)  ^{-2}\right\}  ^{1/\left(
1-2p\tau_{1}\right)  },
\]
Finally, replacing $c$ by its expression, yields%
\[
\omega:=\left\{  -\frac{1}{2\tau_{1}^{3}}\frac{\gamma_{1}^{4}}{C_{2}%
^{2}\left(  D_{1}C_{1}\right)  ^{2p\tau_{1}}}\left(  \int_{0}^{1}%
s^{-1/p+1}K^{2}\left(  s\right)  ds\right)  \left(  \int_{0}^{1}s^{-\tau_{1}%
}K\left(  s\right)  ds\right)  ^{-2}\right\}  ^{1/\left(  1-2p\tau_{1}\right)
},
\]
as sought.

\end{document}